

\magnification\magstep1

\def\makeatletter{\catcode`\@11 }
\def\makeatother{\catcode`\@12 }

\input mssymb

\makeatletter
\font\teneuf=eufm10
\font\seveneuf=eufm7
\font\fiveeuf=eufm5
\newfam\euffam
\textfont\euffam=\teneuf
\scriptfont\euffam=\seveneuf
\scriptscriptfont\euffam=\fiveeuf
\def\frak{\relaxnext@\ifmmode\let\next\frak@\else
 \def\next{\Err@{Use \string\frak\space only in math mode}}\fi\next}
\def\goth{\relaxnext@\ifmmode\let\next\frak@\else
 \def\next{\Err@{Use \string\goth\space only in math mode}}\fi\next}
\def\frak@#1{{\frak@@{#1}}}
\def\frak@@#1{\noaccents@\fam\euffam#1}
\def\relaxnext@{\let\next\relax}
\def\noaccents@{\def\accentclass@{0}}
\makeatother

\font\footnotefont=cmr7
\newcount\footnotecount
\def\myfootnote#1{\global\advance\footnotecount1 
	\footnote{$^{\number\footnotecount}$}%
		 {\footnotefont\openup-3pt{#1}\vskip-0.5cm}}

\font\bigbold=cmbx12 

\line{\bigbold Martin Goldstern\myfootnote{Supported by DFG grant Ko
490/7-1, and by 
	the Edmund Landau Center for research in Mathematical
	Analysis, supported by the Minerva Foundation (Germany)}
				\hfill \rm November, 1992/January, 1993}
\smallskip
\leftline{\bigbold Miroslav Repick\'y\hfill\rm Version of August 1993}
\smallskip
\leftline{\bigbold Saharon Shelah$^{1,}$\myfootnote{Publication 487.}
}
\smallskip
\leftline{\bigbold Otmar Spinas\myfootnote{Supported by the Basic Reasearch
        Foundation of the Israel Academy of Sciences 
        and the Schweizer Nationalfonds}
}
\bigskip
{\bigbold\centerline{ON TREE IDEALS} }
\bigskip

{\leftskip0.12\hsize\rightskip0.12\hsize

	\noindent {\bf Abstract.}\ \ Let $l^0$ and
	$m^0$ be the ideals associated with Laver and
	Miller forcing, respectively. We show that
	${\bf add }(l^0) < {\bf cov}(l^0)$ and
	${\bf add }(m^0) < {\bf cov}(m^0)$ are
	consistent.  We also show
	that both Laver and Miller forcing collapse
	the continuum to a cardinal $\le {\frak
	h}$.

}

\baselineskip1.14\baselineskip

\def \a {\alpha } \def \g {\gamma } \def \l {\langle } \def \r
{\rangle }

\bigskip \bigskip

\newcount\thmcount
\def\prop#1#2:{\advance\thmcount1 \edef#1{\the\thmcount}\medskip
		 \noindent{\bf #1. #2:}}
\def\fct{{^\omega \omega}}

\def\begindent{\par\smallskip\begingroup\advance\parindent by 1cm\relax}
\def\endent{\smallskip\endgroup}
\def\extend{{}^{\frown}}
\def\limpl{\rightarrow}
\def\on{{\restriction}}
\def \c {{\frak c}} \def \o {\omega }
\def\extend{{}^{\frown}}
\def\limpl{\rightarrow}

\def\ite #1 {\item{(#1)}}
\def\cut{\cap}
\def\Sk{{\Bbb S}}
\def\Ml{{\Bbb M}}
\def\Lv{{\Bbb L}}
\def\newline{\hfil\break}

{\bf Introduction and Notation}\ In this paper we investigate the
ideals connected with the classical tree forcings introduced by Laver [La]
and Miller [Mi]. Laver forcing $\Lv$ is the set of all trees $p$ on ${^{<
\o } \o }$ such that $p$ has a stem and whenever $s\in p$ extends
$stem(p)$ then $Succ_p(s):=\{ n: s \, \hat {\, } \, n \in p \} $ is
infinite. Miller forcing $\Ml$ is the set of all trees $p$ on ${^{<
\o } \o }$ such that $p$ has a stem and for every $s\in p$ there is
$t\in p$ extending $s$ such that $Succ_p(t) $ is infinite. The set of
all these splitting nodes in $p$ we denote by $Split(p)$. For any $t\in
Split(p)$, $Split_p(t)$ is the set of all minimal (with respect to
extension) members of $Split(p)$ which properly extend $t$. For both
$\Lv$ and  $\Ml$ the order is inclusion.

The Laver ideal $\ell ^0$ is the set of all $X\subseteq {^\o \o }$ with
the property that for every $p\in \Lv$ there is $q\in \Lv$ extending $p$
such that $X\cap [q]=\emptyset $. Here $[q]$ denotes the set of all
branches of $q$. The Miller ideal $m^0$ is defined analogously, using
conditions in $\Ml$ instead of $\Lv$. By a fusion argument one easily
shows that $\ell ^0$ and $m^0$ are $\sigma $-ideals.

The  additivity (${\bf add}$) of any ideal is defined as the minimal
cardinality of a family of sets belonging to the ideal whose union
does not. The covering number (${\bf cov}$) is defined as the least
cardinality of a family of sets from the ideal whose union is the
whole set on which the ideal is defined -- ${^\o \o }$ in our case.
Clearly $\o _1 \leq $ $ {\bf add}(\ell ^0)\leq $ ${\bf cov}(\ell ^0) \leq \c $ and $\o _1
\leq $ ${\bf add}(m ^0)\leq $ ${\bf cov}(m^0) \leq \c $ hold.

The main result in this paper says that
there is a model 
 of ZFC where ${\bf add}(\ell ^0)<$ ${\bf cov}(\ell ^0)$ and ${\bf add}(m
^0)<$ ${\bf cov}(m ^0)$ hold. The motivation was that by a result of Plewik
[Pl] it was known that the additivity and the covering number of the
ideal connected with Mathias forcing  are the same and they are
equal to the
cardinal invariant ${\frak h}$ -- the least cardinality of a family of
maximal antichains of ${\cal P}(\o )/{\it fin} $ without a common
refinement. On the other hand, in [JuMiSh] it was shown that
${\bf add}(s^0)<$ ${\bf cov}(s^0)$ is consistent, where $s^0$ is
Marczewski's ideal
-- the ideal connected with Sacks forcing $\Sk$. Intuitively, $\Lv$ and
$\Ml$ sit somewhere between Mathias forcing and $\Sk$. In [GoJoSp] it was
shown that under Martin's axiom ${\bf add} (\ell ^0)$ $=$ ${\bf add}(m^0)=\c $,
whereas this is false for $s^0$ (see [JuMiSh]).

The method of proof for ${\bf add}(s^0)<$ ${\bf cov}(s^0)$ in [JuMiSh] is the
following: For a forcing $P$
denote by $\kappa (P)$ the least cardinal to which forcing with $P$
collapses the continuum. In [JuMiSh] it is shown that ${\bf add}(s^0)\leq
\kappa (\Sk)$. In [BaLa] it was shown that in $V^{\Sk_{\o _2}}$ $\kappa
(\Sk) =\o _1$ holds -- where
$\Sk_{\o _2}$ is the countable support iteration of length $\o _2$ of
$\Sk $. Hence $V^{\Sk_{\o _2}}\models
$ ${\bf add}(s^0)=\o _1 $. On the other hand, a L\"owenheim-Skolem argument
shows that $V^{\Sk_{\o _2}}\models $ ${\bf cov}(s^0)= \o _2$.

Our method of proof is similar. Denoting by $P_{\omega_2}$
a  countable support iteration of length $\o _2$ of $\Lv$ and  $\Ml$
(each occuring on a stationary set),
 in $\S 2$ we prove the following:

\bigskip

{\bf Theorem}
$$V^{P_{\o _2}}\models \o _1 = {\bf add}(\ell ^0) = {\bf add}(m^0)    <
{\bf cov}(\ell^0)= 
{\bf cov}(m^0)    = \o _2 $$

\bigskip

The crucial steps in the proof are to show that $\kappa (\Lv)$, $\kappa
(\Ml)$ equal $\o _1$ and ${\bf add}(\ell ^0)\leq \kappa (\Lv)$, ${\bf
add}(m^0)\leq
\kappa (\Ml)$ holds.

We will use the standard terminology for set theory and forcing. By
${\frak b}$ we denote the least cardinality of a family of functions in
$^\o \o $ which is unbounded with respect to eventual dominance and
${\frak d}$ will be the least cardinality of a dominating family in $^\o
\o $. Moreover ${\frak p}$ is the least cardinality of a filter base on $([\o
]^\o ,\subseteq ^* )$ without any lower bound, and ${\frak t}$ is the
least cardinality of a decreasing chain in $([\o ]^\o ,\subseteq ^* )$
without any lower bound.  
It is easy to see that $\omega_1 \le {\frak p } \le {\frak t} \le
{\frak b } \le {\frak d} \le {\frak c } $. 

\bigskip \bigskip

\noindent {\bf 1.\ \ Upper and lower bounds}

\bigskip

{\bf 1.1 Theorem} \ {\it (1) ${\frak t}\leq $ ${\bf add}(\ell ^0)\leq $
${\bf cov}(\ell ^0) \leq {\frak b}$

(2) ${\frak p}\leq $ ${\bf add}(m^0)\leq $ ${\bf cov}(m^0)\leq {\frak d}$}

\bigskip

{\it Proof of (1):} \ We have to prove the first and the third
inequality. For the third inequality, let $\langle f_\a :\a < {\frak b} \rangle $ be an unbounded
family. Define
$$X_\a :=\{ f\in {^\o \o }: (\exists ^\infty k)\; f(k)<f_\a (k)\} $$
Clearly $\bigcup \{ X_\a :\a < {\frak b} \} = {^\o \o }$. We claim $X_\a
\in \ell ^0$. Let $p\in \Lv$. We define $q\in \Lv $ as follows:
$stem(q):=$ $stem(p)$, and for any $s$ extending $stem(q)$ we have
$s\in q$ if and only if $s\in p$ and $(\forall k)$ if $|stem(q)|\leq k
<|s|$ then $s(k)\geq f_\a (k)$. Then clearly $q\in \Lv $, $q$ extends $p$ and
$[q]\cap X_\a =\emptyset $.

\def \l {\langle } \def \r {\rangle } \def \seq {{^{<\o }\o }}
\smallskip

In order to prove the first inequality we use the following
notation from [JuMiSh]: Let $Q:=\{ \bar A =\langle A_s : s\in
{^{<\o }\o }\rangle : (\forall s)\;A_s \in [\o ]^\o \} $. For $\bar A\in Q$ we
define a sequence of Laver trees $\l p_s (\bar A) : s\in \seq \r $ as
follows: $p_s (\bar A)$ is the unique Laver tree such that $stem(p_s
(\bar A))=s$ and if $t\in p_s (\bar A)$ extends $s$ then $Succ_{p_s
(\bar A)}(t) =A_t$.

For $\bar A , \bar B \in Q$ we define:
$$\openup6pt\eqalign{
\bar A \subseteq \bar B &\Leftrightarrow (\forall s)\; A_s \subseteq B_s\cr
\bar A \subseteq ^* \bar B &\Leftrightarrow (\forall s)\; A_s
						\subseteq ^* B_s\cr
\bar A \leq ^*\bar B &\Leftrightarrow (\forall s)\; A_s \subseteq ^*
B_s \wedge (\forall ^\infty s )\;A_s \subseteq B_s\cr}$$

Here $\leq ^*$ is a slight but important modification of $\subseteq
^*$ from [JuMiSh].

\bigskip

{\bf Fact 1.2} \ {\it $(Q,\leq ^*)$ is ${\frak t} $-closed.}

\def \t {{\frak t}}
\bigskip

{\it Proof of 1.2}\ Suppose $\l \bar A_\a :\a < \g \r $ where $\g < \t
$ is a decreasing sequence in $(Q,\leq ^*)$. Let $\bar A_\a :=\l A^\a _s :s \in
\seq \r $. Since $\g <\t $ there is $\bar B'= \l B'_s :s\in \seq
\r \in Q$ such that
$(\forall \a <\g )\;\bar B'\subseteq ^* \bar A_\a $. Define $f_\a : \seq
\rightarrow \o $ such that $(\forall s)\;B'_s \setminus f_s(\a )
\subseteq A^\a _s $. Since $\t \leq {\frak b}$ there exists $f:\seq
\rightarrow \o $ such that $(\forall \a )(\forall ^\infty s)\; f_\a
(s)\leq f(s)$. Now let $B_s :=B'_s\setminus f(s)$ and $\bar B :=\l B_s
:s\in \seq \r $. It is easy to check that $(\forall \a <\g )\;\bar B\leq
^* \bar A_\a $.

\bigskip

{\bf Fact 1.3}\ {\it Suppose $X\in \ell ^0$ and $\bar A \in Q$. There
exists $\bar B \in Q$ such that $\bar B\subseteq \bar A$ and $(\forall s
\in \seq )\;[p_s(\bar B)]\cap X=\emptyset $.}

\bigskip

{\it Proof of 1.3:} \ First note that if $D:=\{ p\in \Lv : [p]\cap X
=\emptyset \} $ then $D$ is open dense and even $0$-dense, i.e. for
every $p\in \Lv$ there exists $q\in D$ extending $p$ such that
$stem(q)=$ $stem(p)$. The proof of this is similar to Laver's proof in
[La] that the set of Laver trees deciding a sentence in the language of
forcing with $\Lv $ is $0-$dense: Suppose $p\in \Lv $ has no 
$0-$extension whose
branches are not in $X$. Then inductively we can construct $q\in \Lv $
extending $p$ such that every extension of $q$ has a branch in $X$, 
contradicting $X\in \ell ^0$.
 
Using this it is straightforward to construct
$\bar B$ as desired.

\bigskip

{\bf Fact 1.4:} \ {\it Suppose $X\subseteq {^\o \o }$, $\bar A,\bar B
\in Q$, $ \bar B \leq ^*
\bar A$ and $(\forall s)\;[p_s(\bar A)]\cap X =\emptyset $. Then $(\forall
s)\;[p_s(\bar B)]\cap X =\emptyset $.}

\bigskip

{\it Proof of 1.4:}\ Clearly, if $F\subseteq p_s(\bar B)$ is finite,
then
$$[p_s(\bar B)] =\bigcup \{ [p_t(\bar B)]: t\in p_s(\bar B)\setminus F
\} $$
 But for almost all $t\in p_s(\bar B)$, $p_t(\bar B)$ extends
$p_t(\bar A)$. So clearly $[p_s(\bar B)]\subseteq [p_s (\bar A)]$ and
hence $[p_s(\bar B)]\cap X =\emptyset $.

\def \b {\beta }
\bigskip

{\it End of the proof of 1.1(1).} \ Suppose we are given $\l X_\a :\a
<\g \r $ and $q\in \Lv$, where $\g < \t $ and $(\forall \a )\;X_\a \in
\ell ^0$. Choose $\bar A\in Q$ such that $p_{stem(q)}(\bar A )=q$ and
let $\bar B_0 $ be the $\bar B$ given by 1.3 for $\bar A$ and $X_0$.
If $\l \bar B_\a :\a <\beta \r $ has been constructed for $\b \leq
\g $
and $\b $ is a successor, then choose $\bar B_\b $ as given by 1.3
for $\bar A= \bar B_{\b -1}$ and $X=X_\b $. If $\b $ is a limit, then
by 1.2 choose first $\bar A$ such that $(\forall \a < \b )\;\bar A \leq
^* \bar B_\a $ and then find $\bar B_\b \subseteq \bar A $ as given by
1.3 for $\bar A$ and $X=X_\b $. Finally, if we have constructed $\bar
B_\g =\l B^\g _s : s\in \seq \r $ define $\bar B:= \l B_s :s\in \seq \r
$ by $B_s:= B^\g _s \cap Succ_{q}(s) $ if $s\in q $ extends $stem(q)$
and $B_s:=B^\g _s$ otherwise. It is easy to check that $\bar B\in Q$,
$p_{stem(q)}(\bar B) $ extends $q$ and $(\forall \a < \g
)\;[p_{stem(q)}(\bar B)]\cap X_\a = \emptyset $.

\bigskip

{\it Proof of 1.1(2)} \ The proof is similar to (1). For the third
inequality, let $\l f_\a :\a < {\frak d}\r $ be a dominating family.
Define
$$X_\a :=\{ f\in \fct : (\forall ^\infty k)\;f(k)<f_\a (k)\} $$
\noindent Then $\bigcup \{ X_\a :\a < {\frak d} \} = \fct $ and in an
analogous way as in (1) it can be seen that $X_\a \in m^0$.

In order to prove the first inequality we need the following concept
\relax from [GoJoSp]. Let $R$ be the set of all $\bar P=\l P_s :s\in \seq \r
$ where each $P_s\subseteq \seq $ is infinite, $t\in P_s$ implies
$s\subset t$ and if $t, t' \in P_s$ are distinct then $t(|s|)\ne
t'(|s|)$. Given $\bar P\in R$ we can define $\l p_s(\bar P):s\in \seq
\r $ as follows: $p_s(\bar P)$ is the unique Miller tree with stem $s$
such that if $t\in Split(p_s(\bar P))$ then $Split_{p_s (\bar
P)}(t)=P_t$.

Define the following relations on $R$:

$$\eqalign{\bar P \leq \bar Q &\Leftrightarrow (\forall s )\;p_s(\bar P)
\leq p_s(\bar Q)\cr
\bar P \approx \bar Q& \Leftrightarrow (\forall s )\;P_s =^*
Q_s \wedge (\forall ^\infty s)\; P_s =Q_s\cr
\bar P \leq ^* \bar Q& \Leftrightarrow (\exists \bar P')\;\bar P
\approx \bar P' \wedge \bar P' \leq \bar Q\cr}$$

\bigskip

{\bf Fact 1.5 [GoJoSp, 4.14]}\ {\it Assume $MA_\kappa $($\sigma
$-centered). If $\l \bar P_\a :\a <\kappa \r $ is a $\leq
^*$-decreasing sequence in $R$, then there exists $\bar Q\in R$ such
that $(\forall \a <\kappa )\; \bar Q\leq ^*\bar P _\a $.}

\bigskip

The following two facts have similar proofs as 1.3 and 1.4.

\bigskip

{\bf Fact 1.6} \ {\it Suppose $X\in m^0$ and $\bar P\in R$. There
exists $\bar Q \leq \bar P$ such that $(\forall s)\; [p_s(\bar Q)]\cap X
=\emptyset $.}

\bigskip

{\bf Fact 1.7} \ {\it Suppose $X\in m^0$, $\bar P ,\bar Q \in R$,
$\bar P\leq ^* \bar Q$ and $(\forall s)\;[p_s(\bar Q)]\cap X =\emptyset
$. Then $(\forall s)\;[p_s(\bar P)]\cap X =\emptyset $.}
\tracingpages1
\bigskip

Now using 1.5, 1.6, 1.7 and the well-known result that for all
$\kappa <{\frak p}$ $MA_\kappa $($\sigma $-centered) holds, a similar
construction as in 1.1(1) shows that ${\frak p}\leq $ ${\bf add}(m^0)$.

\bigskip \bigskip

{\bf 2 \ ${\bf add}$ and ${\bf cov}$ are distinct}
\par\nobreak
\bigskip
{\bf Definition 2.1}\ A set $A \subseteq \fct$ is called {\it strongly
dominating} if and only if
$$ (\forall f\in \fct)( \exists \eta\in A)( \forall^\infty k)\;
	f(\eta(k-1))<\eta(k)$$

\bigskip

{\bf Definition 2.2}\  For any set $A \subseteq \fct$, we define
the domination game $D(A)$ as follows:

There are two players, GOOD and BAD.  GOOD plays first.  The game
lasts $\omega$ moves.
$$\vcenter{\offinterlineskip
\ialign{\strut\hfil$#$\hfil&\vrule#&\hfil$#$\hfil\cr
\omit\hfil ~GOOD~\hfil &&\omit\hfil ~BAD~\hfil\cr
\omit  &height 2pt\cr
\noalign{\hrule}
\omit  &height 2pt\cr
s &&\cr
  && n_0 \cr
m_0 && \cr
    && n_1 \cr
m_1 && \cr
\vdots && \vdots \cr
}}$$
The rules are: $s$ is a sequence in ${^{<\omega }}\o $, and the $n_i$
and $m_i$ are natural numbers. (Whoever breaks these rules first,
loses immediately).

The GOOD player wins if and only if
\begindent
\itemitem{(a)}For all $i$, $m_i>n_i$.

\itemitem{(b)}The sequence $s\extend m_0\!\extend m_1\!\extend \cdots$
is in $A$.

\endent

\bigskip

{\bf Lemma 2.3} \ {\it Let $A \subseteq \fct $ be a Borel set. Then
the following are equivalent:

\itemitem{(1)}There exists a Laver tree $p$ such that $[p] \subseteq A$.

\itemitem{(2)}$A$ is strongly dominating.

\itemitem{(3)}$\hbox{GOOD}$  has a winning strategy in the game $D(A)$. }

\bigskip

{\it Remark:}  \ Strongly dominating is not the same as
dominating.  For example, the closed set
$$A:= \{ \eta\in \fct: (\forall k)\; \eta(2k)=\eta(2k+1) \}$$
is dominating but is not strongly dominating.

\bigskip

{\it Proof of 2.3}\
We consider the following condition:
\begindent
\item{(4)}(For all $F:{^{<\omega}}\o \times \omega \to \omega )(\exists
\eta\in A)( \forall^\infty k)( \forall i\le k )\; \eta(k) > F(\eta\on k,
i)$.

\endent

We will show (1) $\limpl $ (2) $\limpl$ (4) $\limpl $ (3) $\limpl $
(1).

\bigskip
(1) $\limpl$ (2) is clear.

\bigskip
(2) $\limpl$ (4): Given $F$, define $f$ by
$$ f(m) := \max \{ F(s,i): i\le m, s\in m^{\le m+1}\} + m$$
$f$ is increasing, $f(m) \ge m$ for all $m$.

Find $\eta$ such that $\forall^\infty k \, \eta(k) > f(\eta(k-1))$.
Then $\eta$ is increasing.  For almost all $k$ we have, letting $m:=
\eta(k-1)$:\newline
 $m \ge k-1$, so $\eta\on k \in m+1^{m+1}$, so by the
definition of $f$ we get  $f(m)\ge F(\eta\on k,i)$ for any $i\le k$.
So $\eta(k) > f(\eta(k-1) \ge F(\eta\on k,i)$.

\bigskip
(4) $\limpl$ (3):  Assume that GOOD has no winning strategy.  Then BAD
has a winning strategy ${\sigma}$ (since the game $D(A)$ is Borel,
hence determined.)

We can find a function $F:{^{<\omega}}\o \times \omega  \to \omega$
such that for all $s, m_0, \ldots, m_k$ we have
$$ {\sigma}(s, m_0, \ldots, m_k) = F(s\extend m_0 \!\extend
\cdots\extend m_k, |s|) $$
Find $\eta\in A$ as in (4).  So  there is $k_0$ such that $\forall k
\ge k_0 \, \eta(k) \ge F(\eta\on k, k_0)$.  So in the  play
$$\vcenter{\offinterlineskip\ialign{\strut\hfil$#$\hfil&\ \vrule#\
							&\hfil$#$\hfil\cr
\omit\hfil ~GOOD~\hfil &&\omit\hfil ~BAD~\hfil\cr
\omit  &height 2pt\cr
\noalign{\hrule}
\omit  &height 2pt\cr
s:=\eta\on k_0 &&\cr
  && n_0 := {\sigma}(s)=F(\eta\on k_0, k_0) \cr
m_0:= \eta(k_0+1) && \cr
    && n_1 := {\sigma}(s, m_0) = F(\eta\on (k_0+1), k_0)  \cr
m_1:=\eta (k_0+2) && \cr
\vdots && \vdots \cr
}}$$
player BAD followed the strategy ${\sigma}$, but player GOOD won, a
contradiction.

\bigskip
(3) $\limpl $ (1):  Let $B$ be the set of all sequences $s\extend
m_0\!\extend m_1\!\extend \cdots$ that can be played when GOOD follows
a specific winning strategy.  Clearly $B \subseteq A$, and for some Laver tree
$p$, $B= [p]$.

\def\o {\omega } \def\a {\alpha } \def\c {{\frak c}} \def \g {\gamma }
\def \b {\beta }

\bigskip

{\bf Lemma 2.4 [Ke]}\ {\it Let $A\subseteq \fct $ be an analytic set.
Then the following are equivalent:

\itemitem{(1)}There exists a Miller tree $p$ such that $[p]\subseteq
A$.

\itemitem{(2)}$A$ is unbounded in $(\fct ,\leq ^*)$.}

\bigskip

{\bf Lemma 2.5}\ {\it (1) Suppose ${\frak b}=\c $. For every dense open
$D\subseteq \Lv$ there exists a maximal antichain $A
\subseteq D$ such that
$$\forall q\in \Lv ([q]\subseteq \bigcup \{ [p]: p\in A \}\
\Rightarrow \ \exists A'\in [A]^{ <\c } \forall p\in A\setminus A' \,
p\perp q ) \eqno(\ast ) $$

(2) The same is true for $\Ml$.}

\bigskip

{\it Proof:} \ Let $\Lv =\{ q_\a : \a <\c\} $.
Inductively we will define a set $S \subseteq \c$ and sequences
$\langle x_{\gamma}:{\gamma}< \c \rangle$ and
$\langle p_{\gamma}:{\gamma}\in S\rangle$. Finally we will let $A= \{ p_\gamma
:\gamma \in S \} $.

Let $0\in S$ and choose $x_0 \in [q_0]$ arbitrarily. 

It can be easily seen that every Laver tree contains $\c $ 
extensions such that every two of them do not contain a common branch.
So clearly we may find $p_0 \in D$ such that $x_0 \not\in [p_0 ]$.

Now suppose that $\l x_\gamma :\gamma < \alpha \r $ and $\l p_\gamma :
\gamma \in S\cap \a \r $ have been constructed for $\a < \c $.

First choose $x_\a \in [q_\a ]$ arbitrarily, but such that, if 
$[q_{\a }] \not\subseteq \bigcup \{ [p_\gamma ]:\gamma < \a \} $ then 
$x_\a \not\in
\bigcup \{ [p_\gamma ]:\gamma < \a \} $.

In order to decide whether $\a \in S$ or not we distinguish the following two 
cases:

\noindent {\bf Case 1:}  $q_{\a }$ is compatible
with some $p_\gamma $, ${\gamma} <\a $. In this case ${\a }\notin S$.

\noindent {\bf Case 2:} $q_{\a }$ is incompatible with all
$p_\gamma $, ${\gamma} < \a $.  Now we let ${\a }\in S$, and we define 
$p_\a $ as follows:

By Lemma 2.3 for each $ \g  \in
\a $ we may find $f_\g :  \o  \rightarrow \o $ such that

$$(\forall \eta \in [p_\g ]\cap [q_\a ])( \exists ^\infty k)\; \eta (k)
\leq f_\g (\eta (k-1)) \eqno(\ast \ast )$$

By our assumption on ${\frak b}$ there exists a strictly increasing $f$
which dominates all
the $f_\g $'s. Now define $p_{\a } '\in \Lv$ as follows:
  $stem(p_\a ')$ $=$
$stem(q_\a )$, and for $t\in p_{\a } '$,
if $t \supseteq stem (p_\a ')$ and $|t|=:n$
we require 

$$ Succ_{p_{\a } '}(t) = Succ_{q_{\a } }(t) \cut
	[f(t(n-1)), \infty)$$
\noindent Clearly
$p_\a '\in \Lv$, $p_\a ' \subseteq q_\a $, and by $(\ast \ast )$ and our
assumption on $f$ we conclude
$[p_\g ]\cap [p_\a '] =\emptyset $ for every $\g < \a $.

By the remark above that every Laver tree contains $\c $ 
extensions such that every two of them do not contain a common branch, 
we may find $p_\a \in D$ such that
$p_\a $ extends $p_\a '$ and $[p_\a ]$ and $\{ x_\g :\g \leq
\a \} $ are disjoint.

This finishes the construction. Now let
$A := \{p_{\gamma }:{\gamma} \in S\}$. 

Since every $q_{\a }$ is either
compatible with some $p_\g $, $\g < {\a }$ (in case 1)
or contains the condition $p_{\a }$ (in case 2) and for $\a \ne \g $ with
$\a ,\g \in S$ we have $[p_\a ]\cap [p_\g ]=\emptyset $ we conclude that 
$A$ is a maximal
antichain.  

$A$ also satisfies condition $(\ast )$: Let
$q=q_{\a }$. By construction, if $[q_\a ]\not\subseteq \bigcup \{ [p_\g ] : 
\g \in S\cap \a \} $ then $[q_\a ]\not\subseteq \bigcup \{ [p_\g ] : 
\g \in S \} $. 

\smallskip

The proof of (2) is analogous, but instead of Lemma 2.3 we use 2.4.

\def \k {\kappa }

\bigskip

{\bf Lemma 2.6} \ {\it  Suppose ${\frak b}={\frak c}$.
	Then ${\bf add}(\ell^0)\leq \k (\Lv)$ 
and ${\bf add}(m^0)\leq \k (\Ml)$.}
%

\bigskip

{\it Proof:}\ We may assume $\k (\Lv)< \c $. Let $\dot f$ be a
$\Lv$-name such that $\Vdash _{\Lv} \, ``\dot f:\k (\Lv )\rightarrow \c
\hbox{ is onto}$''. For $\a < \k (\Lv)$ let
$$D_\a :=\{ p\in \Lv : (\exists {\beta} )\; p\Vdash _{\Lv} \, \dot
f(\a )=\b \} $$
For $p\in D_\alpha$  will write ${\beta}_p = {\beta}_p(\alpha)$ for
the unique ${\beta}$ satisfying $p \Vdash _{\Lv} \, \dot
f(\a )=\b  $.

 Clearly $D_\a $ is dense and open. So we may choose a
maximal antichain $A_\a \subseteq D_\a $ as in Lemma 2.5. Let
$$X_\a :=  {^\o \o }\setminus \bigcup \{ [p] : p\in A_\a \} $$
Then $X_\a \in \ell ^0$. We claim that $X=\bigcup_{\a <\k
(\Lv)} X_\a \not\in \ell ^0$. Suppose on the contrary $X\in \ell ^0$.
So we may find $q\in \Lv$ such that $[q]\cap X=\emptyset $ and hence
$[q]\subseteq \bigcup \{ [p] : p\in A_\a \} $ for each $\a $. By the
choice of $A_\a $ each of the sets
$$ B_\alpha := \{ {\beta}_p(\alpha): p\in A_\alpha, \hbox{ $p$ compatible
with $q$}\}$$
is bounded in $\c$. Since $\c $ is regular by our assumption
${\frak b}=\c $ we can find $\nu < \c$  such that for all $\alpha <
\kappa(\Lv)$, $B_\alpha \subseteq \nu$. So easily conclude that
$$q \Vdash _{\Lv} \, \hbox{``range$(\dot f) \subseteq \nu < \c $''}$$
This is a contradiction.

The proof for $\Ml$ is similar.

\bigskip

{\bf Theorem 2.7}\  $\kappa (\Lv) \leq {\frak h}$ {\it and} $\kappa (\Ml)
\leq {\frak h}$.

\def \h {\frak h}
\bigskip

{\it Proof:}\ We prove it only for $\Lv$. The proof for $\Ml$ is very
similar.  We work in $V$. Let $\langle {\cal A}_\a :\a < {\bf
h}\rangle $ be a family of maximal almost disjoint families such that,
\begindent
\ite 1  if $\a <\b <\c$ then ${\cal A}_\b $ refines ${\cal A}_\a $
\ite 2  there exists no maximal almost disjoint family refining all
	the ${\cal A}_\a $
\ite 3  $\bigcup \{ {\cal A}_\a : \a < \h \} $ is dense in
	$([\o ]^\o,\subseteq ^*)$
\endent
 That such a sequence exists was shown in [BaPeSi].

Since $\h $ is regular, for every $p\in \Lv$ there exists $\a <\h $
such that for each $s\in Split(p)$  there is $A\in {\cal
A}_\a $ with $A\subseteq ^* Succ_p(s)$. Hence, writing $\Lv_\a $ for
the set of those $p\in \Lv$ for which $\a $ has the property just
stated, we conclude $\Lv=\bigcup \{ \Lv_\a :\a < \h \} $.

For each $A\in {\cal A}_\a $ choose ${\cal B}_A=\{ B^A(p):p\in \Lv\}$,
an almost disjoint family on $A$.

Now we will define $\Lv_\a ':=\{ q^\a(p):p\in \Lv_\alpha \} $ such that
$q^\a (p) $
extends $p$ for every $p\in \Lv_\alpha $ and $p_1\not=p_2$ implies
$q^\a (p_1)\perp q^\alpha(p_2)$.
For $p\in \Lv_\alpha$,  $q^\a(p)$ will be defined as follows:

{\parindent=0cm\leftskip=0.1\hsize\rightskip=0.1\hsize
 For each $s\in Split( p) $ let
$C^\a _{s}(p):= Succ_{p} (s)\cap B^A(p) $ where $A\in {\cal A}_\a $ is
such that $A\subseteq ^* Succ_{p }(s) $.
So clearly $C^\a _{s}(p)$ is infinite. Now
$q^\a (p) $ is the unique Laver tree $ \subseteq p$ satisfying
 $stem(q^\a (p))=stem(p)$ and for each $s\in Split(q^\a (p))$
we have $Succ_{q^\a (p) }(s)=C^\a _{s }(p)$.

}
It
is not difficult to see that $\Lv_\a '$ has the stated properties.

Now we are ready to define a $\Lv$-name $\dot f$ such that $\Vdash
_{\Lv}\, ``\dot f :{\h } ^V \rightarrow \c ^V $ is onto'':
For each $p\in \Lv_\alpha$, let
$\{ r^\a _{\xi }(p): \xi < \c \} \subseteq \Lv$ be a maximal antichain
below $q^\a(p)$, and
define $\dot f$ in such a way that $r^\a _{\xi }(p)\Vdash _{\Lv} \,
``\dot f (\a )=\xi $''.  As $\bigcup \{ \Lv_\a ' :\a <\h \} $ is dense
in $\Lv$, it is easy to check that $\dot f$ is as desired.

\bigskip

{\bf Theorem 2.8} \ {\it Let $\omega_2 = S_{\Ml}\dot\cup S_{\Lv}$,
where the sets $S_{\Ml}$ and $ S_{\Lv}$ are disjoint and stationary. 
Let $(P_\alpha, Q_\alpha: \alpha < \omega_2)$  be a
countable support iteration of length $\o _2$  such that for all
$\alpha$ we have $\Vdash_{P_\alpha} Q_\alpha = \Ml$ whenever
$\alpha\in S_{\Ml}$, and 
$\Vdash_{P_\alpha} Q_\alpha = \Lv$ otherwise. 
Also suppose that $V$ satisfies CH. Then in $V^P$, ${\frak h}=\o _1$
holds.}

\bigskip

{\it Proof:} \ Both $\Ml$ and $\Lv$ have the property $(*)_1$ of
[JuSh].  (For $\Lv$, this was proved in [JuSh] and for $\Ml$ this was
proved in [BaJuSh].)  [JuSh] also showed that this property is preserved
under countable support iterations, so also $P_{\omega_2}$ has this
property.  Hence the reals of $V$ do not have measure zero in $V^P$,
so from $\h  \leq {\frak s} \leq {\bf unif}({\cal L}) $ (where ${\frak
s}$ is the  splitting number and {\bf unif}$({\cal L})$ is the
smallest cardinality of a set of reals which is not null) we get the
desired conclusion. 
%
%

\bigskip

{\bf Theorem 2.9} \ {\it Let $P_{\omega_2}$ be as in 2.8.
Then $$V^{P_{\omega_2}} \models 
 \o _1 = {\bf add}(\ell ^0) = {\bf add}(m^0) < 
	{\bf cov}(\ell^0)= 
	 {\bf cov}(m^0)    = \o _2 $$}

\bigskip

{\it Proof:} \ Since $\Lv$  adds a dominating real, we have
$V^{P_{\omega_2}} \models {\frak b } = {\frak c} $, so by 2.6, 2.7 and
2.8, 
it suffices to prove that the covering
coefficients are $\o _2$ in the respective models. The proof of this is
 similar to the proof of [JuMiSh, Thm1.2] that ${\bf cov}$ of the
Marczewski ideal is $\o _2$ in the iterated Sacks' forcing model.

We give the proof only for $\ell^0 $. Suppose $\langle X_\a :\a <\o _1
\rangle \in  V^{P _{\o _2}}$ is a sequence of $\ell ^0$-sets. In
$V^{P _{\o _2}}$ let $f_\a :\Lv \rightarrow \Lv $ be such that for
every $p\in \Lv $, $f_\a (p)$ extends $p$ and $[f_\a (p)]\cap X_\a
=\emptyset $. Since $P _{\o _2}$ has the $\o _2$-chain condition, by a
L\"owenheim-Skolem argument it is possible to find $\g <\o _2$ such
that
$$\langle f_\a \on \Lv ^{V_{\gamma} } : \a <\o _1 \rangle \in V^{P _\g }$$
where $V_{\gamma}:= V^{P_{\gamma}}$. 
Moreover, it is possible to find such a ${\gamma}$ in $S_{\Lv}$. 
 We
claim that the  Laver real $x_\g $ (which is added by $Q_{\gamma} =
\Lv^{V_{\gamma}}$)  is not in $\bigcup _{\a <\o
_1}X_\a $, which will finish the proof. Otherwise, for some $p\in
\Lv_{\g \o _2}$ where $\Lv_{\g \o _2}:= \Lv_{\o _2}/G_\g $ and some $\a
<\o _1$ we would have: $p\Vdash  \, x_\g \in X_\a $. But letting $q:=p(\g
)\in \Lv $ and letting $r(\g ):=f_\a (q)$ and $r(\b ):=p(\b )$ for $\b >
\g $ we see that $r\Vdash \, x_\g \not\in X_\a $, a contradiction.

\vfill\bigskip\bigskip\bigskip

\centerline{{\bf References}}

\advance\rightskip1cm minus 1cm

\bigskip

\item{[BaJuSh]}T.~Bartoszynski, H.~Judah, S.~Shelah, Cicho\'n's
Diagram, to appear in the Journal of Symbolic Logic.

\item{[BaPeSi]}B. Balcar, J. Pelant, P. Simon, The space of
ultrafilters on $N$ covered by nowhere dense sets, Fund. Math.,
110(1980), 11-24.

\item{[BaLa]}J.E. Baumgartner and R. Laver, Iterated perfect set
forcing, Ann. Math. Logic, 17(1979), 271-288.

\item{[GoJoSp]}M. Goldstern, M. Johnson and O. Spinas, Towers on
trees, Proc. AMS, to appear.

\item{[JuMiSh]}H. Judah, A. Miller, S. Shelah, Sacks forcing, Laver
forcing and Martin's axiom, Arch. Math. Logic, 31(1992), 145-161.

\item{[JuSh]}H. Judah and S. Shelah, The Kunen-Miller chart, J. Symb.
Logic, 55(1990), 909-927.

\item{[Ke]}A. Kechris, A notion of smallness for subsets of the Baire
space, Trans. AMS, 229(1977), 191-207.

\item{[La]}R. Laver, On the consistency of Borel's conjecture, Acta
Math., 137(1976), 151-169.

\item{[Mi]}A. Miller, Rational perfect set forcing, Contemporary
Mathematics, vol.31(1984), edited by J.E. Baumgartner, D. Martin and 
S. Shelah, 143-159.

\item{[Pl]}S. Plewik, On completely Ramsey sets, Fund. Math.
127(1986), 127-132.

\bigskip\bigskip\bigskip

{\bf Addresses:} \parindent2cm\def\\#1{\item{{#1}}}

\\{\it Martin Goldstern} 2.~Mathematisches Institut,
Freie Universit\"at Berlin, Arnimallee 3, \newline
14195 Berlin, Germany.
\newline
e-mail: {\tt goldstrn@math.fu-berlin.de}

\smallskip

\\{\it Miroslav Repick\'y} Matematick\'y \'ustav SAV, Jesenn\'a 5, 04154
Ko\v sice, Slovakia
\newline
e-mail: {\tt repicky@ccsun.tuke.cs}
\smallskip

\\{\it Saharon Shelah} Institute of Mathematics, Hebrew University of
Jerusalem, \newline Jerusalem, Israel.
\newline
e-mail: {\tt shelah@math.huji.ac.il}
\smallskip

\\{\it Otmar Spinas} Departement Mathematik, ETH-Zentrum, 8092 Z\"urich, 
Switzerland {\it and} \newline
Institute of Mathematics, Hebrew University of Jerusalem,\newline
Jerusalem, Israel (current address).
\newline
e-mail: {\tt spinas@math.ethz.ch}
\bye